\documentstyle[12pt,amssymb]{article}
\newtheorem{Theo}{Theorem}
\newcommand{\www}{{available on www.maths.ox.ac.uk/$\sim$ zilber}}

\newcommand{\Q}{{\mathbb Q}}
\newcommand{\dd}{{\partial}}
\newcommand{\size}{{\rm size}}
\newcommand{\f}{{\rm f}}
\newcommand{\K}{{\bf K}}
\newcommand{\M}{{\bf M}}
\newcommand{\Z}{{\mathbb Z}}
\newcommand{\R}{{\mathbb R}}

\newcommand{\krn}{{\rm  ker}}
\newcommand{\rem}{{\bf  Remark }}
\newcommand{\EC}{{\cal  EC}}
\newcommand{\N}{{\mathbb N}}
\newcommand{\A}{{\cal A}}

\newcommand{\std}{{\rm st}}
\newcommand{\chr}{{\rm char }}
\newcommand{\trd}{{\rm tr.d.}}
\newcommand{\ex}{{\rm ex}}
\newcommand{\ld}{{\rm lin.d.}}
\newcommand{\HH}{{\cal H}}
\newcommand{\card}{{\rm card\ }}

\newcommand{\C}{{\mathbb C}}
\newcommand{\ssn}{\section}
\newcommand{\bt}{\begin{Theo}}
\newcommand{\et}{\end{Theo}}
\newcommand{\inv}{^{-1}}
\newcommand{\lb}{\label}
\newcommand{\be}{\begin{equation}}
\newcommand{\ee}{\end{equation} }
\newcommand{\ra}{\rangle}
\newcommand{\la}{\langle}
\newcommand{\subs}{\subseteq}
\newcommand{\subsf}{\subseteq_{fin}}

\begin{document}
\title{ Analytic and pseudo-analytic structures}
\author{ B. Zilber\\ University of Oxford }  \date{2 April 2001}

\maketitle

One of the questions frequently asked nowadays about  
model theory is
whether it is still logic. The reason for asking the question is
mainly that  more and more of model theoretic research focuses on
concrete  mathematical fields, uses extensively their tools and attacks
their inner problems. Nevertheless the logical roots in the case of model theoretic
geometric stability theory  are not only clear but also remain very
important in all its applications.

This line of research started with the notion of a $\kappa$-categorical first
order theory,
which quite soon mutated into the more algebraic and less logical notion
of a $\kappa$-categorical structure.\\

A structure $\M$ in a first order  language $L$ is said to be {\bf categorical in cardinality $\kappa$}
if there is  exactly one, up to isomorphism, structure  of cardinality
$\kappa$ satisfying the $L$-theory  of  $\M.$

In other words, if we add to ${\rm Th}(\M)$ the (non first-order) statement that
the cardinality of the domain of the structure is $\kappa,$ the description becomes
categorical.\\

The principal breakthrough, in the mid-sixties, from which  stability
theory started was the answer to J.Los' problem\\

{\bf The  Morley Theorem} {\em A countable theory which is categorical in
one uncountable cardinality is categorical in all uncountable cardinalities.}
\\

{The basic examples of uncountably categorical structures in a
countable language are:}\\

(1) Trivial structures (the language allows  only equality);\\

(2) Abelian divisible torsion-free groups; Abelian groups of prime
    exponent (the language allows $+, =$); Vector spaces over a
    (countable) division ring\\

(3) Algebraically closed fields in language $(+,\cdot, =)$ .\\

Also, any structure definable in one of the above is uncountably categorical
in the language which witnesses the interpretation.

The structures definable in algebraically closed fields, for example,
are effectively objects of algebraic geometry. \\

As a matter of fact  the main logical problem after answering the
 question of J.Los was 
{\em what properties of $\M$ make it $\kappa$-categorical for
 uncountable $\kappa?$}\\

{ The answer is now reasonably clear:}

{\em The key factor is  measurability by a dimension and high
homogeneity of the structure.}\\ 

This gave rise to (Geometric) Stability Theory,  the theory studying 
structures with
good dimensional and geometric properties (see [Bu] and [P]). When
applied to fields, the
stability theoretic approach in many respects is very close to Algebraic Geometry.\\

The abstract dimension notion for finite $X\subset \M$ mentioned above could be best  understood
by examples:   \\

(1a) Trivial structures: \ \ {\bf size of } $X$; \\

(2a) Abelian divisible torsion-free groups; Abelian groups of prime
    exponent; Vector spaces over a division ring:  \  {\bf linear dimension of } $X$;\\

(3a) Algebraically closed fields: \  {\bf transcendence degree}
 $\trd(X).$ 
 
 Dually, one can classically define another type of dimension using
 the initial one: 
 $$\dim V=\max\{ \trd (\bar x) \ | \  \bar x\in V\}$$ for
 $V\subs \M^n$ an algebraic variety. The latter type of dimension notion
is called in model theory the {\bf Morley rank}.
\\

The last example can serve also as a good illustration of the
significance of homogeneity of the structures. So, in general, the
transcendence degree makes good sense in any field, and there is 
quite a reasonable dimension theory for algebraic varieties over a
field. But the dimension theory in arbitrary fields fails if we want
to consider it for wider classes of definable subsets, e.g. the images
of varieties under algebraic mappings. In algebraically closed fields
any definable subset is a boolean combination of varieties, by
elimination of quantifiers, which eventually is the consequence of the
fact that algebraically closed fields are existentially closed 
in the class of fields. The latter effectively means high homogeneity,
as an existentially closed structure absorbs any amalgam with another
member of the class.      \\

One of the achievements of stability theory is the establishing of
some hierarchy of types of structures that allows to say which
ones are more 'analysable' (see [Sh]).\\

The next natural question to ask is {\em whether
there are 'very good' stable structures which are not reducible to
(1) - (3) above?}\\ 

The initial hope of the present author in [Z1],  that any uncountably categorical
structure comes from the classical context (the trichotomy conjecture), was based on the general
belief that logically perfect structures could not be overlooked in
the natural progression of mathematics. Allowing some philosophical
licence 
here, this was also a
belief in a strong logical predetermination of basic mathematical structures.

As a matter of fact this turned out to be true in many
cases. Specifically for 
{\em Zariski geometries}, which are defined as the structures with a
good dimension theory and nice topological properties, similar to the
Zariski topology on algebraic varieties
(see  [HZ]).

Another situation where this principle works, is the context of
o-minimal structures (see [PS]).
\\

Powerful applications of the result on Zariski geometries and of the underlying methodology
were found by Hrushovski [H3],[H4]. This not only lead to  new
and independent solutions to 
some Diophantine problems,
Manin-Mumford and Mordell-Lang (the functional case) conjectures, but
also  a new geometric vision of these.\\

Yet the trichotomy conjecture proved to be false in general as
Hrushovski  found a source of a great variety of 
counterexamples.

We analyse below the Hrushovski construction, purporting to
answer the question of whether the counterexamples it provides
dramatically overhaul the trichotomy conjecture or if there is a way to
save at least the spirit of it. As the reader will find below the
author is inclined to stick to the second alternative.\\ \\

\ssn{  Hrushovski construction of new structures } 

The main steps:\\

Suppose we  have a, usually elementary, class of structures $\HH$
with a good dimension notion $d(X)$ for finite subsets of the structures.
We want to introduce a new function or relation on $\M\in \HH$ so that the new
structure gets a good dimension notion.
  
The main principle, which Hrushovski found will allow us to do this, is
that of the free fusion. That is, the new function should be
related to the old structure in as a free way as possible. At the same
time we want the structure to be homogeneous. He then
found an effective way of writing down the condition: {\em the number of
 explicit dependencies in $X$ in the new structure must not be greater
than the size (the cardinality) of $X.$} 

The explicit $L$-dependencies
on $X$ can
be counted as $L$-codimension, $\size(X)-d(X).$  The explicit dependencies
coming with a new relation or function are the ones given by
simplest 'equations', basic formulas. 

So, for example, if we want a
{\em new unary function $f$ on a field} (implicit in [H2]), the condition should be 

\be  \lb{e1} \trd(X\cup f(X))-\size(X)\ge 0,\ee
since in the set $Y=X\cup f(X)$ the number of  explicit
field dependencies is $\size(Y)-\trd(Y),$ and the number of explicit dependencies in terms of
$f$ is $\size(X).$

If we want, e.g., to put a {\em new ternary relation $R$ on a field}, then
the condition would be  
 \be  \lb{e2} \trd(X)- r(x)\ge 0,\ee
where $r(X)$ is the number of triples in $X$ satisfying $R.$

The very first of Hrushovski's examples (see [H1]) introduces just a
{\em new
structure of a ternary relation}, which effectively means putting new
relation on the trivial structure. So then we have
\be  \lb{e3} \size(X)-r(X)\ge 0.\ee

If we similarly introduce an automorphism $\sigma$ on the field ({\em
difference fields,} [CH]), then we have to count

\be  \lb{e4} \trd(X\cup \sigma(X))-\trd(X)\ge 0,\ee

and the inequality here always holds.

Similarly for {\em
differential fields} with the differentiation operator $D$ (see [Ma]), where we always have

\be  \lb{e5} \trd(X\cup D(X))-\trd(X)\ge 0.\ee

The left hand side in each of the inequalities (\ref{e1}) -
(\ref{e5}), denote it
$\delta(X),$ is a counting 
function, which is called {\bf predimension},  as it satisfies some of 
the basic properties of the dimension notion.\\

At this point we have carried out the first step of the Hrushovski
construction, that is:\\

(Dim) we  introduced the class $\HH_{\delta}$ of the
structures with a new function or relation, and the extra condition

$$ {\rm (GS)} \ \ \ \ \ \ \ \ \ \ \delta(X)\ge 0 \mbox{ \ \ \ for all finite } X.$$

(GS) here stands for 'Generalised Schanuel', the reason for which will
be given below.
The condition (GS) allows us to introduce another counting function with 
respect to a given structure $\M\in \HH_{\delta}$ 

$$ \dd_M(X)=\min\{ \delta(Y): X\subs Y\subsf M\}.$$

We also need to adjust the notion of embedding in the class for further 
purposes. This is the {\bf strong embedding}, $\M\le {\bf L},$  meaning that
$\dd_{M}(X)=\dd_L(X)$ for every $X\subsf \M.$\\ 

The next step is 

(EC) Using the inductiveness of the class construct an existentially
closed structure in $(\HH_{\delta}, \le).$

If the class has the amalgamation property, then the existentially closed 
structures are sufficiently homogeneous. Also {\em $\dd_M(X)$ for 
existentially closed $\M$  becomes a dimension notion}. 

So, if also the class $\EC$ of existentially closed structures  is axiomatisable, one can rather easily check
that the 
existentially closed 
structures are $\omega$-stable. This is the case for examples
(\ref{e1}) - (\ref{e3}) and (\ref{e5})
above.

In more general situations the e.c. structures may be unstable, but still
with a reasonably good model-theoretic properties.

Notice that though condition (GS) is trivial in examples (\ref{e4}) -
(\ref{e5}), the derived dimension notion $\dd$ is
non-trivial. In both examples $\dd(x)>0$ iff the corresponding rank of
$x$ is 
infinite (which is  the SU-rank
 in algebraically closed difference  fields 
 and  the Morley rank, in differentially closed fields).\\

Notice that the dimension notion $\dd$ for finite subsets,
similarly to the example~(3a), gives rise to
a dual dimension notion for definable subsets $S\subs \M^n$ over a
finite set of parameters $C:$
$$\dim(S)=\max \{ \dd(\{x_1,\dots,x_n\}/C): \ \la x_1,\dots,x_n\ra \in S\}.$$

(Mu) This stage originally had been considered prior to (EC), but as one 
easily sees, it can be equivalently introduced after.

We want to find now a finite Morley rank structure as a substructure
(may be non-elementary) of a structure $\M\in \EC.$ 
In fact an existentially
closed $\M$  would be of finite Morley rank, if '$\dim(S)=0$' is
equivalent to '$S$ is finite'. But in general $\dim(S)$ may be zero for
some infinite definable subsets $S,$ e.g. the set $S=\{ x\in \M: f(x)=0\}$
is one such  in example (\ref{e1}) : 'some equations have too many solutions'.   

To eliminate the redundant solutions Hrushovski
introduces a counting function $\mu$ for the maximal allowed size of
potentially Morley rank 0 subsets. Then
$\HH_{\delta, \mu}$  is the subclass of structures of $\HH_{\delta}$ satisfying the bounds
given by $\mu.$
 Equivalently, since existentially closed structures are
universal for structures of $\HH_{\delta},$ so 
$\HH_{\delta, \mu}$ is the class  
of substructures of existentially closed structures $\M,$
satisfying the bound by $\mu.$  

Inside this class we
can just as well carry out  the construction of existentially closed structures
$\M_{\mu}.$ Again, if the subclass has the amalgamation property
and is first order definable, then an existentially closed
substructure 
$M_{\mu}$ of this subclass is of finite Morley
rank, in fact strongly minimal in cases (\ref{e1}) - (\ref{e3}). 
It is also important for the further discussion that $\M_{\mu}\subs \M.$\\

The infinite dimensional structures emerging after step (EC) in natural classes we call
{\em natural Hrushovski structures}. Some but not all of them lead
after step (Mu) to finite Morley rank structures.

It follows immediately from the construction, that the class of
natural Hrushovski structures is singled out in $\HH$ by three
properties: the generalised Schanuel property (GS), the  property of existentially closedness
(EC) and the property (ID), stating the
existence of $n$-dimensional subsets for all $n.$ 

It takes a bit more model theoretic analysis, as is done in [H1],
to prove that in examples (\ref{e1})-(\ref{e3}), and in many others,
(GS), (EC) and (ID) form a complete set of axioms.  \\ \\

 Since Hrushovski found the counterexamples, the main question that
 has arisen is
 whether the pathological structures demonstrate the failure of
 the general principle or if there is a classical context which explains
 the counterexamples.

We  now want to try and find grounds for the latter.
 \\

We start with one more example of Hrushovski construction.
\ssn{\ Pseudo-exponentiation} 
Suppose we want to put a new function
$\ex$ on a field $K$ of characteristic zero, so that $\ex$  is a 
homomorphism between the additive and the multiplicative groups of the
field: $$\ex(x_1+x_2)=\ex(x_1)\cdot \ex(x_2).$$
Then the corresponding predimension on new structures 
$\K_{\ex}=(K,+,\cdot,\ex)$ must be
$$\delta(X)=\trd(X\cup \ex(X))-\ld(X)\ge 0, \ \ \ \ \ {\rm (GS)}$$    
where  $\ld (X)$ is the linear dimension of
$X$ over $\Q.$

Equivalently (GS) can be stated  
$$\mbox{assuming that $X$ is linearly independent over
$\Q,$ \ \ \  } \trd(X\cup \ex(X))\ge \size(X),$$
which in case $\K$ is the field of complex numbers and $\ex=\exp$ is
known as the {\em Schanuel conjecture} (see [La]).\\

Start now with the class $\HH(\ex/\std)$ consisting of structures
$\K_{\ex}$ satisfying (GS) and with the additional property that the
kernel $\krn=\{ x\in K: \ex(x)=1\}$ is a cyclic subgroup of the
additive group of the field $K,$ which we call a {\bf standard kernel} (see [Z2]). This class is non-empty and can be
described as a subclass of an elementary class defined by omitting
countably many types.

The subclass $\EC(\ex/\std)$ of existentially closed substructures of 
$\HH(\ex/\std)$ is
first order axiomatisable inside $\HH(\ex/\std)$ (see [Z3] for the
proof  about a similar structure).  

By the obvious analogy with the structure $\C_{\exp}=(\C, +,\cdot,
\exp)$ on the complex numbers 
we conjecture that $\C_{\exp}$ is one of the structures in $\EC({\ex}/\std).$
And we want to find the condition that might
single out the isomorphism type of $\C_{\exp}$ among the other structures
in the class. We do this by introducing an extra step to the
Hrushovski construction, when it is possible, which comes after (EC). 
This is applicable in a wide variety of classes:\\

(P) Consider a
$\dd$-independent set $C$ of cardinality $\kappa$ (embeddable in an existentially closed structure)
and let $E(C)$ be
a structure {\bf prime over $C$ in class $\EC$}. \\

When a prime structure over $\dd$-independent $C,$ $\card C=\kappa,$ exists
 and is unique over $C,$ we call $E(C)$ {\bf the $\kappa$-canonical
structure.}\\

The nice thing about having $\kappa$-canonical structures is that we
get a link with the logical background again: the notion is a good
analogue of  $\kappa$-categoricity, or rather the strong minimality
concepts (see also [Le] for recent developments in this direction).\\   

We prove

{\em If the $\kappa$-canonical structure $E(C)$ exists and the language is
countable, we have   $\card E(C)=\kappa$ and  $E(C)$ is
$\aleph_0$-quasi-minimal, i.e. any definable subset of $E(C)$ is either
countable, or the complement of a countable. Moreover, when
$\kappa>\aleph_0,$ for a definable
$S\subs E(C)^n$ \ $S$ is countable \ iff \ $\dim(S)=0.$}    \\

It follows from a well-known theorem of Shelah (see [Sh] and [Bu]), that  canonical
structures exist if the class $\EC$ is first order axiomatisable, complete
and $\omega$-stable, e.g. in examples
(\ref{e1}) - (\ref{e3}) and (\ref{e5}). There is no prime structure in
the class of algebraically closed difference fields
(example~(\ref{e4})) due to the fact that the theory is not complete.

In spite of the fact that $\EC(\ex/\std)$ is not axiomatisable and even
interprets the ring of integers we managed to prove in [Z2] that

{\em There exists a weaker (non unique) version of $\kappa$-canonical structure for any infinite cardinal $\kappa$ in
$\HH(\ex/\std)$ which is $\aleph_0$-quasi-minimal.  }\footnote{ {\em Added
in proof.} We now proved that there is an
$L_{\omega_1,\omega}(Q)$-sentence axiomatizing a subclass of
$\HH(\ex/\std)$ with unique model in every infinite cardinal, and the
models are $\aleph_0$-quasi-minimal.
We don't
know whether the models are canonical in the above sense.}

This finally brings us to \\ 

{\bf Conjecture}  {\em $\C_{\exp}$
is  the canonical  structure of cardinality $2^{\aleph_0}$ in
$\HH(\ex/\std)$ if we put $\ex$ to be $\exp$}.\\

The conjecture (and even the one stating that $\C_{\exp}\in \EC(\ex/\std)$)
is obviously  stronger than the Schanuel conjecture.
Another  consequence of the conjecture is the fact that $\C_{\exp}$ is
existentially closed. We show in [Z2] that this is equivalent to the
statement\\

EC(exp): \ {\em Any non obviously contradictory system of equations over $\C$ in terms of
$+,$ $\cdot$ and $\exp$  has a solution in $\C.$}\\

 The definition of {\bf non obviously contradictory system} is quite
 similar to ones (implicitly) formulated for other classes in the form of axiom
 schemes, e.g. ACFA(iii) in [CH] for algebraically closed difference
 fields, and see also examples in section~\ref{s3}. 
 We are not going to give the definition  here,
 but a good example of such a system is an equation of the form
 $t(x)=0,$ where $t(x)$ is a term in $+,  \cdot, \exp$ over $\C$ which is not
 of the form $\exp(s(x))$ for some other term $s(x).$ Such an equation by
 the conjecture should have a solution in $\C.$ I have learned, while
 writing this paper, that such was exactly a conjecture by S.Schanuel
which was proved by W.Henson and L.Rubel using Nevanlinna theory (see
 [HR]). A.Wilkie gave a rather simple proof of the solvability of an equation
 of the form $\sum_{i\le N}q_i(x)e^{p_i(x)}=0,$ for $q_i$ and $p_i$
 polynomials in one variable, $N>1,$ based on the rate of growth argument.
The general case is
 open, but some research suggesting the truthfulness of the conjecture
 can be traced in the literature. See also the discussion below.

Notice also that one can easily replace the exponentiation by other
classical functions and observe similar effects, including
corresponding versions of the Schanuel conjecture. 

Based on the analysis of pseudo-exponentiation one would like to
conclude hypothetically that \\

1. Basic Hrushovski structures have  analytic prototypes.\\

2. The statement of the Schanuel conjecture along with its analogs is
an 
intrinsic
   property of classical analytic functions, probably responsible for
   a good dimension theory of the corresponding structures on the complex
   numbers. \\

3. Another basic property of classical analytic structures is the
EC-property: {\em Any non-obviously contradictory system of equations
over the structure  has a solution.}\\

Of course, we don't possess  technical means to check the truthfulness
of the conjectures. But the general picture drawn above in view
of the conjectures can be tested in simpliest examples.    
\\

\ssn{   Analytic interpretations}\lb{s3} 

{\bf New ternary relation}\\

Let $g: \C^3 \to \C$ be an
entire function with the properties:\\

(GS[$R$])  If a system of $n+1$ 
equations of
the form $g(v_{i_k},v_{j_k},v_{l_k})=0$ $(k=1,\dots,n+1),$ with
$v_{i_k},v_{j_k},v_{l_k}\in \{ v_1, \dots,v_n\},$
has a 
solution $\la a_1,\dots, a_n\ra$ in $\C,$ then at least two of the triples  
$\la a_{i_k},a_{j_k},a_{l_k}\ra$ coincide. \\

(EC[$R$]) Let $\{ \la x_{j},y_{j},z_{j}\ra : j\le n+m\}$ be  distinct
triples, where each of $x_{j}, y_j,z_j$ is either a complex constant or one of
variables $v_1, \dots, v_n,$ but not all three of them constants.     
Then the system  $$\{ g(x_{i},y_{i},z_{i})=0: i\le n\}\cup  
 \{ g(x_{i},y_{i},z_{i})\neq 0: n< i\le n+m\}$$  of $n$ 
equations and $m$ inequalities has a solution in $\C,$ provided no
$k$ of the $n$ equations have  less than $k$ explicit variables. \\  

It is easy to see that
if we interpret the ternary relation 
 $R(v_1,v_2,v_3)$ as $g(v_1,v_2,v_3)=0,$ 
by GS[$R$] we get condition (\ref{e3}) for $\delta,$  and by
EC[$R$]  the existential closedness is satisfied. Then ID[$R$] follows
from the fact that $\dd$-closure of a finite set is countable, since
$n$ independent analytic equations in the complex $n$-space have only
countably many non-singular solutions. \\

{\bf Problem 1}. (i) Construct an entire function $g$  such that
GS[$R$] and EC[$R$] hold.  

(ii) Prove that $(\C,R)$ is canonical in this case, i.e. the structure
is prime over its basis (of cardinality $2^{\aleph_0}).$ 
\\

\rem One can construct $g$ satisfying (GS[$R$]), using an argument of A.Wilkie, as
$$g(v_1,v_2,v_3)=\sum_{i_1,i_2,i_3\in \N}a_{i_1,i_2,i_3}v_1^{i_1}v_2^{i_2}v_3^{i_3}$$
with complex coefficients $a_{i_1,i_2,i_3}$ algebraically independent
and very rapidly decreasing.\\

{\bf  New functions on a field.} \\

This is another class  mentioned above. One can easily write down in a
first order way the following two schemes of axioms:\\

(GS[\f]) Let $V\subs K^{2n}$ be  a variety over $\Q$ in variables
$x_1,\dots,x_n,y_1,\dots,y_n.$ If $\dim V<n,$ then there is no  
point $\la x_1,\dots,x_n,y_1,\dots,y_n\ra$ in $V$ with $y_i=\f(x_i)$
and $x_i\neq x_j$
for all distinct $i,j\in \{ 1,\dots,n\}.$

(EC[\f])  Let $V\subs K^{2n}$ be  an irreducible variety over $K$ in variables
$x_1,\dots,x_n,y_1,\dots,y_n$ such that 

(i) $V$ is not contained in a hyperplane given by an equation of the
form $x_i=x_j$ for $i<j\le n,$ or $x_i=c,$ for $c\in K;$

(ii) for any $0<i_1<\dots<i_k\le n$ the dimension of
$V_{i_1,\dots,i_k},$ the projection of $V$ onto
$(x_{i_1},\dots,x_{i_k},y_{i_1},\dots,y_{i_k})$-space, is not less
than $k.$

Then there is a   
point $\la x_1,\dots,x_n,y_1,\dots,y_n\ra$ in $V$ with $y_i=\f(x_i),$
for all $i\in \{ 1,\dots,n\}.$\\

It has been proved in [Z6] that\\

{\em  {\rm GS[\f]} and {\rm EC[\f]} along with {\rm ID[\f]}, the axiom of
infinite $\dd$-dimensionality, determine a complete $\omega$-stable theory
(of Morley rank $\omega$). The theory has canonical model $\K_{\f}$ in every
cardinality $\kappa.$}\\

{\bf Problem 2} (i) Construct an entire holomorphic function $\f:\C\to
\C$ satisfying GS[\f] and EC[\f] \ for $K=\C.$

(ii) Prove that
$\C_{\f}$ is isomorphic to  canonical $\K_{\f}$ of cardinality
$2^{\aleph_0}.$  \\ \\

Notice that we get basically the same theory, with minor changes, if we
change GS[\f] and, correspondingly, the counting function $\delta$ in
(\ref{e1}) to:\\

 (GS'[\f]) \ $\trd(X\cup \f(X))-\size(X)\ge 0,$ provided  $X$
does not contain certain elements, say, $0.$   \\
 
A.Wilkie in [W]  proves  that an entire analytic function given as
$$\f(x)=\sum_{i\ge 0} \frac{x^i}{a_i}$$
 with $a_i$ very
rapidly increasing integers,
satisfies (GS'[\f]) and  (EC[\f]), if in the latter $V$ is defined over $\Q .$
\\ \\

{\bf New differentiable functions on a field}.\\

To get better understanding of the 'analytic' nature of previous examples
we consider the class $\HH(F)$ of fields $K$ of characteristic zero with a 
collection $F=\{ \f^{(i)}: i\in \Z\}$ of unary functions on $K.$

We first introduce, 
for a finite $I\subs \Z$
$$\delta_I(X)=\trd(X\cup \bigcup_{i\in I} f^{(i)}(X))- \size( X)\cdot
\size( I),$$ 
{\bf the $I$-predimension} of $X\subsf K.$\\

Then the predimension of $X\subsf K$ is

$$\delta(X)=\min\{ \delta_I( X):I\subsf \Z\}.$$ 

We then, as usual, introduce the subclass satisfying first order
definable condition\\

\hspace{2cm} (GS[$F$])\ \ \ \ \ \ \ \ \ \  \ \ \ \ \ $\delta(X)\ge 0,$\\
 
and after going through the construction stage (EC) find out that the
resulting structures satisfy first order axiom scheme EC[$F$], similar to EC[\f].\\

If we add the corresponding ID[$F$], {\em the first order theory
defined by 
 $${\rm GS}[F]\cup {\rm EC}[F]\cup{\rm ID}[F]$$
 is $\omega$-stable of
rank $\omega,$ and the reduct of the theory to $(+,\cdot, \f),$ where
$\f=\f^{(0)} ,$ is just the theory in the previous example.}\\ 

We then consider for each $i\in \Z$ a definable function
$$g^{(i)}(x_1,x_2)= \left\{ \begin{array}{l} 
\frac{\f^{(i)}(x_1)-\f^{(i)}(x_2) }{x_1-x_2},
\mbox{ if } x_1\neq x_2\\ \\
\f^{(i+1)}(x_1),\mbox{ \ \ \ \  if } x_1=x_2\end{array}\right.$$ 

We want to introduce on $\K_F,$ a model of the theory, a Zariski type
topology $\tau.$ Consider a formal 'completion' \ $\bar K=K\cup \{ \infty\},$ \
which can be viewed as the projective line over the field $K.$ Define
{\bf basic $\tau$-closed subsets of $\bar K^n$} to be

(i) all Zariski-closed subsets  of   $\bar K^n$ for every $n\in \N,$

(ii) for $n=2$ the 'closures' of graphs of $\f^{(i)}:$

$$\{ \la x,y\ra \in \bar K^2: (x\in K \wedge y=\f^{(i)}(x))\vee x=\infty\}$$
 
and

(iii) for $n=3$ the 'closures' of graphs of $g^{(i)}:$

$$\{ \la x_1,x_2,y\ra \in \bar K^3: (x_1,x_2\in K \wedge
y=g^{(i)}(x_1,x_2))\vee 
x_1=\infty\ \vee x_2=\infty \}.$$

Notice that (ii) and (iii) would actually be the closures of the
graphs in the complex topology, if $K=\C$ and $f^{(i)}$ were
holomorphic functions.

Define now the family of {\bf $\tau$-closed sets } to be the minimal
family of  subsets of $\bar K^n,$ all $n,$ containing the
$\tau$-closed subsets and closed under intersections,
finite unions and projections $\bar K^{n+1}\to \bar K^n.$
 We prove in [Z6] that\\

{\em The topology $\tau$ is compact in the  sense that 
any  family with the finite intersection property of $\tau$-closed sets has a non-empty
intersection, and the 
projection of a closed set is closed (but not Hausdorff).} 

 Notice that there is
no DCC for closed sets.\\

We can then view $K$ as a locally compact space. By definition, 
$f^{(i)}$ and $g^{(i)}$ are continuous functions on $K,$ in the sense
of the
topology, and $f^{(i+1)}$  satisfies the definition of the
derivative of $f^{(i)}.$ Moreover, one can then carry out complex
style analysis  and 'non-standard analysis' on $K,$ as shown in [Z4]. 

In particular,  using the notion of {\em infinitesimal elements} in
$K^*\succ K,$ in the case $\chr K=0,$
one gets the Taylor formula for any definable function $h$: {\em given  $x\in K,$  $n\in \N$ and an infinitesimal
$\alpha,$ there is an infinitesimal $\beta$ such that
$$h(x+\alpha)=\sum_{0\le k\le
n}\frac{h^{(k)}(x)}{k!}\cdot\alpha^k+\alpha^n\beta,$$} 
where $h^{(k)}$ are the derivatives of $h.$
In case of positive characteristic $p$ it holds only for $n<p.$\\

\rem It is worth mentioning that in the example of {\em new function on a
field} one could consider $\f$ to be a one-to-one function with the same
predimension formula (\ref{e1}), which is in fact
the original Hrushovski  example in [H2] in the case of an equi-characteristic
 pair of fields. In this case  {\em we can not} expand the
structure to include derivatives under a reasonable topology as
above, and this agrees very well with the fact that there is no
transcendental analytic one-to-one function. But any one-to-one Hrushovski structure, by the way of
construction, is a substructure of a natural Hrushovski structure above.\\ 

\ssn{ Mixed characteristics structures}

There are examples of Hrushovski structures which have no analytic
proptotype. Such an example is {\em the fusion of two
fields}:

Let $\HH$ be the class of two-sorted structures $(L,R)$ with both
sorts
fields, $\chr L=p,$ $\chr R=q,$ in the language of fields for both
sorts and an extra binary relation between the two sorts, interpreted
as a bijective
function $\f: L\to R.$    

One then introduces the predimension function for $X\subsf L$
\be \lb{e6} \delta(X)= \trd(X)+\trd(\f(X))- \size(X)\ge 0\ee
 
After going through steps (Dim) and (EC) one gets an $\omega$-stable
theory of rank $\omega$ which is called the {\bf $\omega$-stable
fusion of two fields.} Notice that the pull-back of the field
structure of $R$ to $L$ gives a second field strcuture on $L,$ which
was  the original Hrushovski example. More importantly, Hrushovski
showed that after carrying out step (Mu) one gets a strongly minimal
fusion of two fields.  \\


Let  $\HH$ again be the class of two-sorted structures with both sorts
fields, $\chr L=0,$ $\chr R=p,$ with the condition,  that $L$ is a
valued field with the residue field $R$ and the valuation group $\Z.$
So  there are also definable mappings $v: L\to \Z$ and $\rho: L_0\to
R,$ where $v$ is the valuation, $L_0=\{ x\in L: v(x)\ge 0\},$ the
valuation ring, and $\rho$ the residue ring homomorphism. 

Following the freeness principle one comes to the following
predimension function, analogous to (\ref{e4}) and (\ref{e5}), for the
class:

Let $v(p)=1,$ and for $X\subsf L$ \ \
$[X]=\{ p^zx: x\in X, z\in \Z\ \mbox{ and } v(x)=-z\}.$ Then 
$$\delta(X)= \trd(X)- \trd(\rho[X]).$$
Like in (\ref{e4}) and (\ref{e5}), $\delta(X)\ge 0$ holds automatically in
this class, so $\HH_{\delta}=\HH$ and every embedding is strong.

Then, completing step (EC), one gets the subclass of existentially
closed structures, which are elementary equivalent to {\em maximal unramified
extensions of the $p$-adic field $\Q_p.$} A more complex structure,
with a definable automorphism $\sigma,$ is studied in [BM] and is
identified there as {\em the field of Witt vectors with an automorphism}.\\

Another  class of mixed charcteristics structures, {\em algebraically
closed valued fields}, with the stability
flavor is studied in [HHM], and there is a strong evidence that this
example is of the pattern under consideration.\\

We present here a series of examples of mixed characteristics of a
different kind.\\

Let $A$ be a 1-dimensional torus or an elliptic curve defined over a
field $F_0$ of
characteristic $p,$ by which we mean just the group scheme, so for
every field $F$ containing $F_0$ there is an algebraic group $\A(F)$
(written multiplicatively) of
$F$-points 
of $\A.$ Let $\HH(\A)$ be a class of two sorted structures $(D,A),$
where $A=\A(F)$ is as above, and with the language of Zariski closed
relations on $A,$ $D$ is  a field in
characteristic zero considered with family of relations and operations
$\Sigma,$ all 
definable by polynomial equations (with parameters) in the field and, at least,
containing the additive operation, so $D$ is a group with an extra
structure.  

Notice, that both sorts are then strongly minimal structures and so have
corresponding dimension notions $\dim_D$ and $\dim_A.$     

Let also the language of the class to contain a function symbol $e,$ which
is interpreted as a surjective homomorphism  $$e:D\to A.$$

We can do the Hrushovski construction for this class putting for
$X\subsf D$
$$\delta(X)= \dim_D(X)+\dim_A(e(X))-\ld(X),$$
where $\ld(X)$ is the dimension in the additive group, i.e. the linear
dimension over $\Q.$

It is easy to see that the resulting structure will have $D$ and $A$
existentially closed, so $F$ is an algebraically closed field. Also, the
kernel $\krn = \{ x\in D: e(x)=1\}$ is an additive subgroup of $D,$
and by algebraic reasons $\krn$ is elementarily equivalent to:

\hspace{0.6cm}  $\Z[1/p],$ \ \ in case $A$ is the multiplicative
group of the field;

or to

\hspace{0.6cm} $\Z[1/p]^{2-i}\times \Z^{i},$ 
in case $A$ is an elliptic curve with
Hasse invariant  $i\in \{0,1\}.$ \\

Here   $\Z[1/p]$ is 
the additive subgroup of rationals with
denominators $p^n,$ $n\in \N,$
if $p>0,$ and is just the group of integers, if $p=0.$ \\

(See also about elliptic curves in [Ha] and the characterisation of elementary equivalence for abelian groups
in [EF]).\\

It is well known that in the given language $A(F)$ is biinterpretable
with the algebraically closed field $F,$ perhaps expanded by finitely
many constants. Thus for $Y\subsf A$ we have \ $\dim_A(Y)=\trd(Y/C),$ with  
the transcendence degree on the right calculated by means of the biinterpretation,
and $C$ is the set of constants.

Because of trichotomy results for definability in algebraically closed
fields (see [R]), following also from later work [HZ], there are essentially
two possibilities to consider for $D:$

\hspace{0.6cm}(i) $D$ is an algebraically closed field;

\hspace{0.6cm}(ii) $D$ is a vector space over a field $K$ of characteristic zero.\\

Thus in the case (i) the inequality for $\delta$ takes form
\be \lb{e6.1} \trd(X)+\trd(e(X))- \ld(X)+d\ge 0,\ee
where $d$ is a non-negative integer depending on the constants needed
for the interpretation.

In case (ii) we have
\be \lb{e7} \ld_K(X)+\trd(e(X))-\ld(X)+d\ge 0,\ee
where $\ld_K(X)$ is the dimension in the sense of the $K$-vector space.\\

The existentially closed structures in the classes bear rather close
similarity to {\em universal coverings of corresponding algebraic varieties}
$A,$
and the corresponding kernel $\ker$ plays the role of {\em the
fundamental group of the variety}, denoted usually $\pi_1(A).$ Both model
theoretically and algebraic geometrically the following two types of the kernel are 
especially interesting: 

\hspace{0.6cm} the case of minimal group satisfying the
condition above, i.e. precisely the group $$\Z[1/p]\mbox{ \ \ or \ \ }
\Z[1/p]^{2-i}\times \Z^i,$$ depending on $A,$ which we call {\bf
standard kernel};

\hspace{0.6cm} the case of  minimal {\em algebraically compact group}
(see [EF]) satisfying the
condition above, i.e. precisely the group $$\prod_{l\neq p \mbox
{ prime } } \Z_l\mbox
{ \ \ \ or \ \ \ }
\prod_{l\neq p \mbox{ prime }} \Z_l^{2-i}\times \prod_{l \mbox{ prime
}} \Z_l^i,$$ depending on $A,$ which we call {\bf
compact kernel}. (Here $\Z_l$ stands for the additive group of
$l$-adic numbers).

The standard kernel corresponds to the, so-called, {\em  topological}
$\pi_1(A),$
and the compact kernel corresponds to the, so-called, {\em  algebraic}
$\pi_1(A)$ [M].

So even in the mixed characteristic case the construction results in
something which, even if  not analytic, has a definite analytic
flavour. But anyway, we hardly know what is the right generalisation
of analyticity in many cases and there is hardly a
satisfactory theory in the case of non-Archimedean valued fields.\\

The case corresponding to formula~(\ref{e6.1}) with $p=0$ is
quite similar to the case of the field with
pseudoexponentiation.

The case   corresponding to formula~(\ref{e7})  is
studied in [Z3] mainly under assumption that the characteristic of $F$
is zero. Notice that the corresponding
structure is quite rich in expressive power, in particular, in sort
$A$ we can 'raise to powers $k\in K$' by putting for $x\in A$
$$x^k=e(k\ln(x)),$$
where $\ln(x)$ is an $e$-pull-back of $x$ in $D,$ so this is a
'multivalued operation'.   We call $\EC$-structures corresponding to
this case {\em the group $A$ with 'raising to $K$-powers'}.  It is also
interesting that, given $k_1,\dots,k_n\in K,$
\be \lb{H} H=k_1\krn +\dots +k_n\krn\mbox{ \ and \ }e(H)=e(k_1\krn) \cdot\dots
\cdot e(k_n\krn)\ee
are definable subgroups of $D$ and $A,$ correspondingly, and when the
kernel is standard these are finite group-rank subgroups. If $n=1$
and $k_1=N\inv$, where $N\in {\cal Z}$ is a non-standard integer, the
group $$e(\frac{1}{N}\krn)$$ is a torsion subgroup  of $A.$

We can prove the existence of the canonical
structure   in all the classes, including pseudo-exponentiation, with
fixed compact kernel (see [Z2]). 
But it remains open\\

{\bf Problem 3} Prove the existence of  $\kappa$-canonical structures
for classes above with standard kernel.\\

Notice that when $K=\Q,$ formula~(\ref{e7}) with $d=0$ holds for any
structure in the class, and  the class $\EC$ (with all possible kernels) is axiomatisable, complete and
superstable. We call the structures in the class {\bf group covers
of $\A.$}
Even for this class Problem 3 is non trivial and has been proved for
the case $\A$ is the multiplicative group of a field of characteristic zero very
recently in [Z7] using some rather subtle field arithmetic results
and some model theoretic techniques of Shelah's. 

We believe that the answer
for group covers with standard kernel is crucial for understanding the
general case. \\ \\

{\bf Further model theoretic properties and a Diophantine conjecture}.   \\

Any deeper model theory of most of the above mentioned classes depends on
the following conjecture about intersections in semi-Abelian varieties
(charcteristic zero case):\\

{\em Given a semi-Abelian variety $B$ over $\Q$  and an algebraic
subvariety $W\subs B$ over $\Q$ 
there is a finite
collection $\tau(W,B)$ of proper semi-Abelian subvarieties of $B$ such
that: for any algebraic subgroup $S\le B$ and any irreducible component 
$U$ of $S\cap W$

either \ \ $\dim U=\dim W+\dim S-\dim B,$ \ \ or \ \ $U\subs T$ for
some $T\in \tau(W,B).$ }\\ \\

The conjecture, as a matter of fact, is of  Diophantine type:

{\em The conjecture on intersections in semi-Abelian
varieties implies the Mordell-Lang (and Manin-Mumford) conjecture for
number fields.}

This can be proved {\em ab initio} as in [Z5], but a nicer, model
theoretical way is to see directly as in [Z3]
that \\

{\em under the assumption that the conjecture is true, 
the structure $A$ with raising to $K$-powers is superstable, 
and the definable finitely generated or torsion subgroups of $A$ are
locally modular.}\\

In fact, the proof of this theorem goes via elimination of
quantifiers to the level of existential formulas ({\em near model
compactness,} as typical for Hrushovski structures), which yields local modularity of the kernel. Hence one
gets local modularity of subgroups in (\ref{H}), which by the
definition of local modularity implies the Mordell-Lang
statement. It remains to notice that any finitely generated or torsion subgroup
of $A$ is embeddable in a group $e(H)$ of (\ref{H}) for a right choice of $K$
and $k_1,\dots, k_n.$ \\

The conjecture effects the fields with pseudo-exponentiation:\\

 {\em Assume the conjecture on intersection in semi-Abelian
varieties, the case $B=(\C^*)^n$.
Then on any canonical model $\K_{\ex}$ in $\EC(\ex/\std)$ there is a (non-Hausdorff)
locally compact topology under which $\ex$ is continuous and,
moreover, in infinitesimal neighborhoods of any point $a\in \K$ \
$\ex$ can be represented by the Taylor expansion}
$$ \ex(a+x)=\ex(a)\cdot \sum_n\frac{x^n}{n!}.$$

Also under the conjecture the structure of genuine exponentiation (or
rather raising to powers) on the
complex numbers becomes very clear.\\

(See [Z5]) {\em Assume the Schanuel conjecture and the conjecture on
intersections in semi-Abelian varieties (the case
$B=(\C^*)^n$). Then the 
structure  of complex numbers with the (multivalued) 
operations $y=\exp(r\ln x)=x^r$ of raising to real
powers satisfies {\rm [EC]}, i.e. any non-obviously contradictory
system of equations of the form
$$\sum_{r_1,\dots r_n}a_{r_1,\dots r_n}x_1^{r_1}\cdot\dots \cdot x_n^{r_n}=0,$$
with $a_{r_1,\dots r_n}\in \C$ and $r_1,\dots, r_n\in \R,$
has a solution in $\C.$  }\\

The proof is based on a theory developed by D.Bernstein,
A.Kushnirenko, B.Kazarnovski and A.Khovanski (see [Kh]).\\

It follows then from the results of preceeding sections that\\

(See [Z3]) {\em Assume the Schanuel conjecture and the conjecture on
intersections in semi-Abelian varieties (the case
$B=(\C^*)^n$). Then the 
structure  of complex numbers with the (multivalued) 
operations of raising to real
powers is superstable.}\\

It is also interesting to notice that the Diophantine conjecture is
equivalent to Hrushovski-Schanuel type inequality~(\ref{e7}) with
$d=0$ and $K={\cal Q},$
a non-standard model of the field of  rationals. In this form we can suggest
the corresponding conjecture for case $p>0:$

\be \lb{e8} \ld_{\cal Q}(X)+\trd(e(X))-\ld_{Q({\cal P})}(X)\ge 0,\ee
where $${\cal P}=\{ p^N: N\in {\cal Z}\}$$ is the subgroup of the multiplicative group of
${\cal Q}$ inducing automorphisms on $A,$ and $Q({\cal P})$ is the subfield of ${\cal Q}$
generated by ${\cal P}.$ 
Thus the corresponding class $\EC$ generalises the class of
algebraically closed difference fields, with commuting generic
automorphisms. So, most probably the theory of
the class is {\em simple}.\\

One can easily rewrite the conjecture in (\ref{e8}) in a more standard
form as follows \\

{\em Given a semi-Abelian variety $A$ over an algebraically closed
field $F$ of characteristic $p$   and an algebraic
subvariety $W\subs A^n$
there are  finite
collections $\lambda(W,A^n)$ of constants  and $\tau(W,A^n)$ of
$n$-tuples $\la f_1 (v_1,\dots,v_m), \dots,f_n (v_1,\dots,v_m)\ra$ of
polynomials with integer coefficients     with the property: 

for any algebraic  subgroup $S\le A^n$ and any irreducible component 
$U$ of $S\cap W$ there are a constant $c\in \lambda(W,A^n),$ an
$n$-tuple $\la f_1, \dots,f_n\ra\in
\tau(W,A^n)$ and
non-negative integers $t_1,\dots,t_m$ such that 
$$S\subs \{ \la x_1,\dots,x_n\ra\in A^n:
x_1^{f_1(p^{t_1},\dots,p^{t_m})}\cdot\dots\cdot x_n^{f_n(p^{t_1},\dots,p^{t_m})}=c\}.$$}

{\bf References}\\

[BM] L.Belair and A.Macintyre, {\em L'automorphisme de Frobenius des
vecteurs de Witt}, C.R.Acad.Sci.Paris, t.331(2000), Serie 1, 1-4\\

[B] S.Buechler {\bf Essential Stability Theory}, Springer, Berlin 1996\\

[CH] Z.Chatzidakis and E.Hrushovski, {\em Model theory of difference
fields}, Preprint 1995\\

[EF] P.Eklof and  E.Fisher, {\em The elementary theory of abelian
groups}, Ann.Math.Logic, 4 (1972), 115-171\\

[Ha] R.Hartshorne, {\bf Algebraic Geometry}, Springer-Verlag, 
Berlin - Heidelberg - New York, 1977 \\

[HR] W.Henson and L.Rubel, {\em Some applications of Nevanlinna theory to
mathematical logic: identities of exponential
functions}, Trans.Amer.Math.Soc., 282 (1984), 1-32. Correction vol.
294 (1986), 381\\

[H1] 
E.Hrushovski, {\em A New Strongly Minimal Set}, Annals of Pure and Applied
Logic, 62 (1993), 147-166 \\

[H2]
------------, {\em Strongly Minimal Expansions of Algebraically Closed Fields},
Israel Journal of Mathematics 79 (1992), 129-151 \\

[H3] 
------------, {\em The Mordell-Lang Conjecture for Function Fields}, Journal
of the AMS, 9 (3) (1996), 667-690 \\

[H4] -----------, {\em Difference fields and the Manin-Mumford
conjecture}, Preprint, 1996\\

[HHM] D.Haskel, E.Hrushovski and D.Macpherson, {\em Independence and
elimination of imaginaries in algebraically closed valued fields},
Preprint, 2000\\

[HZ]  E.Hrushovski, B.Zilber, {\em Zariski  Geometries}. Journal of AMS, 9(1996),
1-56  \\

[Kh] A.Khovanski, {\bf Fewnomials} (In Russian). Fazis, Moscow, 1997\\

[L] S.Lang {\bf Introduction to Transcendental Numbers},
Addison-Wesley, Reading, Massachusets, 1966\\

[Le] O.Lessmann, {\em Dependence relation in non-elementary classes},
Ph.D. thesis, Carnegie Mellon U., 1998\\ 

[M] D.Mumford, {\bf Abelian Varieties,} Tata Instutute, Bombay, 1968\\

[Ma] D.Marker, {\em Model theory of differential fields}, In: {\bf
Model Theory of Fields}, LNL, Springer, 1996\\

[P] A.Pillay, {\bf Geometric Stability Theory}, Clarendon Press,
Oxford, 1996\\

[PS] Y.Peterzil and S.Starchenko, {\em A trichotomy theorem for o-minimal
structures,} Proc.London Math.Soc. (3) 77 (1998), 481-523  \\

[R] E.Rabinovich, {\bf Interpreting a field in a sufficienly rich
incidence system}, QMW Press, 1993\\

[Sh1] S.Shelah, {\em Finite diagrams stable in power}, {\bf Annals of
Math.Logic} 2 (1970), p.69-118\\

[Sh] ---------, {\bf Classification Theory}, revised edition, North-Holland--
Amsterdam--Tokyo 1990\\

[W] A.Wilkie, {\em Liouville functions}, This volume\\

[Z1] B.Zilber, {\em The structure of models of uncountably  categorical  theories}.  In:
{\bf Proc. Int. Congr. of Mathematicians,}  1983,  Warszawa, PWN - North Holland
P.Co., Amsterdam - New York - Oxford, 1984, v.1, 359-368\\

[Z2] ---------, {\em Fields with pseudo-exponentiation}, Submitted,
1999,\\ \www \\

[Z3] ---------, {\em Raising to powers in algebraically closed
fields}, To appear in J.Math.Log. \\

[Z4] ---------, {\em  Quasi-Riemann Surfaces} In: {\bf Logic: from
Foundations to Applications,  
European Logic Colloquium}, eds W.Hodges, M.Hyland, C.Steinhorn and
J.Truss, Clarendon Press, Oxford,  1996, 515-536\\

[Z5] ---------, {\em Exponential sums equations and the Schanuel
conjecture}, J.London Math.Soc.(2) 65 (2002), 27-44\\ 

[Z6] ---------, {\em A theory of a generic function with derivations},
 In:{\bf Logic and Algebra,} 
(ed. Yi Zhang), pp. 85--100, Contemporary Mathematics, Vol 302, AMS, 2002\\

[Z7] ---------, {\em Covers of the multiplicative group of an algebraically closed
field of characteristic zero},  math.AC/0401301

\end{document}